\documentclass{article}
\usepackage{amssymb}
\usepackage{amsmath}
\usepackage{color}

\usepackage[T1]{fontenc}
\usepackage{hyperref}			

\usepackage{amsthm}               
\usepackage{amsmath}             

\newtheorem{teo}{Theorem}[section]
\newtheorem{defi}[teo]{Definition}
\newtheorem{lem}[teo]{Lemma}

\newtheorem{ass}[teo]{Assumption}

\newcommand{\be}{\begin{equation}}
\newcommand{\ee}{\end{equation}}
\newcommand{\ba}{\begin{array}}
\newcommand{\ea}{\end{array}}
\newcommand{\bee}{\begin{eqnarray*}}
\newcommand{\eee}{\end{eqnarray*}}
\newcommand{\bea}{\begin{eqnarray}}
\newcommand{\eea}{\end{eqnarray}}

\newcommand{\comment}[1]{}

\newcommand{\E}{\mathbb E}
\newcounter{algo}[section]
\renewcommand{\thealgo}{\thesection.\arabic{algo}}
\newcommand{\algo}[3]{\refstepcounter{algo}
\begin{center}
\framebox[\textwidth]{
\parbox{0.95\textwidth} {\vspace{\topsep}
{\bf Algorithm \thealgo : #2}\label{#1}\\
\vspace*{-\topsep} \mbox{ }\\
{#3} \vspace{\topsep} }}
\end{center}}
\def\IR{\hbox{\rm I\kern-.2em\hbox{\rm R}}}
\def\IC{\hbox{\rm C\kern-.58em{\raise.53ex\hbox{$\scriptscriptstyle|$}}
    \kern-.55em{\raise.53ex\hbox{$\scriptscriptstyle|$}} }}
    
\newcommand{\RR}{\mathbb{R}}

\def \E {{\mathbb E}}

\usepackage{bbm}				
\def \1{{\mathbbm 1}}

\newcommand{\cS}{{\cal S}}
\newcommand{\cN}{{\cal N}}

\title{Linesearch Newton-CG methods for convex optimization with noise}
\author{
  S. Bellavia\thanks{Department of Industrial Engineering of Florence,
    Universit\`{a} degli Studi di Firenze, Italy. Member of the INdAM Research
    Group GNCS. Email: stefania.bellavia@unifi.it},
  E. Fabrizi\thanks{Department of Mathematics and Computer Science ``Ulisse Dini'',
    Universit\`{a} degli Studi di Firenze, Italy.  {Member of the INdAM Research
    Group GNCS.} Email: eugenio.fabrizi@unifi.it}\hspace*{2pt} and
  B. Morini\thanks{Department of Industrial Engineering of Florence,
    Universit\`{a} degli Studi di Firenze, Italy. Member of the INdAM Research
    Group GNCS. Email: benedetta.morini@unifi.it}
  }

\date{}

\begin {document}
\maketitle

\begin{abstract}
This paper studies the numerical solution of strictly convex unconstrained optimization problems by linesearch Newton-CG methods. We focus on  methods employing inexact evaluations of the objective function and inexact and possibly random gradient and Hessian estimates. 
The derivative estimates  are not required to satisfy
suitable accuracy requirements at each iteration but with sufficiently high probability. Concerning the evaluation of the objective function
we first assume that the noise in the objective function evaluations  is bounded in absolute value. Then, we analyze the case where the  error  satisfies  prescribed dynamic accuracy requirements. 
We provide for both cases  a complexity  analysis and  derive expected iteration complexity bounds.   We finally focus on  the specific case of finite-sum minimization which is typical of machine learning applications.
\end{abstract}

\section{Introduction}

In this paper we consider globally convergent Inexact Newton methods for solving the strictly convex
unconstrained optimization problem
\be\label{problem1}
\min_{x\in\RR^n} f(x). 
\ee
We focus on the Newton method where the linear systems are solved by the Conjugate Gradient (CG) method \cite{HS},
usually denoted as  Newton-CG method,
and on the enhancement of its convergence properties by means of Armijo-type conditions.

The literature on  globally convergent  Newton-CG methods is well established 
as long as the gradient and the Hessian matrix are computed exactly or approximated sufficiently 
accurately in a deterministic way, see e.g., \cite{C, EW, NW}. On the other hand, the research is currently very active for
problems with inexact information on $f$ and its derivatives and possibly such that the 
accuracy cannot be controlled in a deterministic way  \cite{BSV,bg,BGM,BellGuriMoriToin22, belgmt,bkk,bkm,bbn_2017,Berahas,STORM2,noc2,noc1,STORM1,dskkv,frz,LR,PS,mit1,xu2}. 

This work belongs to the recent stream of works and 
addresses the solution of (\ref{problem1}) when the objective function $f$ is computed with noise
and gradient and Hessian estimates are random. Importantly, derivative estimates  are not required to satisfy
suitable accuracy requirements at each iteration but with sufficiently high probability. Concerning the evaluation of $f$
we cover two cases: estimates of $f$ subject to noise that is bounded in absolute value;
estimates of $f$ subject to a controllable error, i.e., computable with a prescribed dynamic accuracy. 
Such a class of problems has been considered in \cite{BellGuriMoriToin22, Berahas}, our contribution consists
in their solution with linesearch Newton-CG method and, to our knowledge, this case
has not been addressed in the literature. We provide two linesearch Newton-CG methods suitable for the class of problems specified above and 
provide bounds on the expected number of
iterations required to reach a desired level of accuracy in the optimality gap.

The paper is organized as follows. In Section \ref{sec2} we give preliminaries on Newton-CG and
on the problems considered. In Section \ref{sec3} we present and  study a linesearch Newton-CG algorithm
where function estimates are subject to a prefixed deterministic noise. In Section \ref{sec4} 
we propose and study a linesearch Newton-CG algorithm where function estimates have controllable accuracy.
In Section \ref{sec5} we consider the specific case where $f$ is a finite-sum which is typical of machine learning applications and compare our approach with  the Inexact Newton methods specially designed for this class of problems given in 
  \cite{bkk,bbn,mit1}.

In the rest of paper  $\|\cdot\|$ denotes the 2-norm. 
{Given symmetric matrices $A$ and $B$,  $A\preceq  B$ means that $B-A$ is 
positive semidefinite. }

\section{Our setting}\label{sec2}
{ In this section we provide  preliminaries on the solution of problem (\ref{problem1}) and the assumptions made.
Our methods belong  to class of the Inexact Newton methods \cite{DES}  combined 
with a linesearch strategy for enhancing convergence properties.
A key feature  is that function, gradient and Hessian 
evaluations are approximated and the errors in such
approximations are either deterministic or stochastic, as specified below.

The Inexact Newton methods considered here are iterative processes where, 
given the current iterate $x_k$, a random approximation $g_k$ to $\nabla f(x_k)$
and a random approximation  $H_k$ to $\nabla^2 f(x_k)$, the trial step $s_k$  satisfies
\be \label{inexact1}
  H_k s_k = -  g_k+ r_k, \quad \|r_k\| \leq \eta_k \| g_k\|,
\ee
for some  $\eta_k\in(0,\bar\eta), 0<\bar\eta<1$, named forcing term.

With $s_k$ and a trial steplength $t_k$ at hand, some suitable sufficient decrease Armijo condition is tested
on $x_k+t_k s_k$. Standard linesearch strategies are applied using the true function $f$. On the other hand, here 
we assume that the evaluation of $f$ is subject to an error. 

If $H_k$ is positive definite we can solve inexactly the linear systems $ H_k s=-  g_k$
using the Conjugate Gradient (CG) method \cite{HS}. The resulting method is denoted as Newton-CG.
If the initial guess for CG is the null vector,  the following properties hold. 
}
\begin{lem} \label{Lemma3a}
{ Suppose that $H_k$ is symmetric positive definite and 
and 
$ s_k $ is the vector in (\ref{inexact1}) obtained by solving $ H_k s=-  g_k$ 
with the CG method and null initial guess.}
{Let $0<\lambda_1\le\lambda_n$ 
such that  $\lambda_1 I \preceq  H_k \preceq \lambda_n I$.}
 Then, there exist constants $\kappa_1, \kappa_2, \beta>0$, such that:
\be \label{hpconstant}
\kappa_1\|g_k\| \leq \| s_k \| \leq \kappa_2 \|g_k\|, \quad -g_k^Ts_k \geq \beta \| s_k \|^2, \qquad \forall k>0,
\ee
which satisfy
\be
\frac{1-\bar \eta}{\lambda_n} \leq \kappa_1 \leq  \kappa_2\leq \frac{1}{\lambda_1}, \quad \beta \geq \lambda_1. \label{kappa12}
\ee
As a consequence
\begin{eqnarray}
&&\beta \kappa_1 \leq \frac{-g_k^Ts_k}{\| s_k \|^2} \frac{\| s_k \|}{\| g_k \|} = \frac{| g_k^Ts_k |}{\| g_k\| \| s_k \|}\leq1 \label{betak1}\\
&&\lambda_n\kappa_1 \geq  1-\bar \eta\\
&&\kappa_1 \lambda_1 \leq \kappa_2 \lambda_1 \leq 1 \label{k11}\\
&&-g_k^Ts_k \geq \beta \kappa_1^2 \| g_k\|^2. \label{betaaa}
\end{eqnarray}

\begin{proof}
{ Lemma 7 in \cite{fg} guarantees that  any step $\hat s$ returned by the CG method applied 
to $H_ks=g_k$ with null initial guess
satisfies $\hat s^TH_k\hat s=-\hat s^Tg_k$. 
}
Then, it holds
\be \label{shs7}
s_k^TH_ks_k = -s_k^Tg_k,
\ee
and by (\ref{convexityH}) we have
\begin{eqnarray} 
-  g_k^T s_k  & =  & s_k^T  H_k s_k \geq \lambda_1 \|s_k\|^2   
\label{bound_prod_a}
\end{eqnarray}
which provides the lower bound on $\beta$ in (\ref{kappa12}). \\
By (\ref{shs7}) it also follows  that
\be \label{slg}
\| s_k \|\leq \lambda_1^{-1}\| g_k \|,
\ee
which provides the upper bound on $\kappa_2$ in (\ref{kappa12}).\\
Moreover, using  (\ref{inexact1})
$$
\| g_k - r_k\| \geq \|  g_k\| - \|r_k\| \geq (1-\eta_k) \|  g_k\| \geq (1-\bar\eta) \|  g_k\|
$$
while  (\ref{inexact1}) and (\ref{convexityH}) gives
$$
\|  g_k - r_k\| = \|  H_k s_k\| \leq \|  H_k \| \| s_k\| \leq \lambda_n \|s_k\|
$$
 and consequently
\be \label{in1}
\|s_k\| \geq \frac{\|  g_k - r_k\|}{\lambda_n} \geq \frac{1-\bar\eta}{\lambda_n} \|  g_k \|,
\ee
which provides the lower bound on $\kappa_1$ in (\ref{kappa12}).\\
Inequalities (\ref{betak1})-(\ref{betaaa}) are direct consequences of (\ref{hpconstant}) and (\ref{kappa12}).
\end{proof}
\end{lem}

\subsection{Assumptions}
We introduce the assumptions on the problem (\ref{problem1}) and on the approximate evaluations of functions, gradients and Hessians.
\begin{ass}  {\rm (smoothness and strong convexity of $f$)}\label{strconv} The function $ f $ is  twice continuously differentiable and there exist some $\lambda_n\ge  \lambda_1 > 0 $ such that the Hessian matrix $\nabla^2 f(x)$ satisfy
\be \label{convexity}
\lambda_1 I \preceq \nabla^2f(x) \preceq \lambda_n I,\;\;\; \forall x \in {\mathbb R}^n. 
\ee
 \end{ass}
\noindent
As a consequence, $f$ is strongly convex with constant $\lambda_1$, i.e.,
$$
f(x)\geq f(y) \nabla f(y)^T(x-y) +\frac{\lambda_1}{2}\| x-y \|^2 \text{ for all } x,y\in\RR^n
$$
and  the gradient of $f$ is Lipschitz-continuous with constant $\lambda_n$, i.e.,
$$
\| \nabla f(x) -\nabla f(y)\| \leq \lambda_n\|x-y \| \text{ for all } x,y\in\RR^n.
$$
Further, letting $ x^* $ be the unique minimizer of the function $f$, for  all $ x \in \mathbb{R}^n $ we have
\be \label{Lemma2}  
\lambda_{1} \|x-x^*\| \leq \|\nabla f(x) \| \leq \lambda_{n} \|x-x^*\|,
\ee
and 
\be \label{conv1}  \frac{\lambda_1}{2 } \|x-x^{*}\|^{2}\leq f(x)-f(x^{*})\leq \frac{1}{2\lambda_1} \|\nabla f(x)\|^{2}. \ee
see \cite[Theorem 2.10]{nesterov}.

As for  approximated evaluations of  the objective function, we consider 
two possible cases. The first one is such that
the value $f(x)$ is approximated with a value $\tilde f(x)$ and the corresponding error is not controllable but 
its upper bound  $\varepsilon_f$ is known.

\begin{ass}\label{ass1f}{\rm (boundness of noise in $f$)}\label{errf} There is a positive scalar $\varepsilon_f$  such that
 \be
 | f(x) - \tilde f(x) |\leq \varepsilon_f, \quad  \forall x\in \RR^n \label{errf}
\ee
\end{ass}

The second case concerns a controllable error between $f(x)$ and $\tilde f(x)$  at the current iteration $x_k$ 
and at the trial iteration $x_k+t_ks_k$.
\begin{ass} \label{ass2f} {\rm (controllable noise in $f$)}\label{efk}
For all $k>0$ and some given positive $\theta$
\begin{eqnarray} \label{errfk}
| f(x_k) - \tilde f(x_k) | &\leq& -\theta t_k s_k^Tg_k,  \nonumber \\
  | f(x_{k}+t_ks_k) - \tilde f(x_{k}+t_ks_k) | &\leq& -\theta t_k s_k^Tg_k. \label{errfk}
	\end{eqnarray}
\end{ass}

The methods we are dealing with are globalized Inexact Newton methods  employing random estimates $g_k$ and $H_k$ of the gradient and the Hessian and noisy values of the objective function. 
Then, they generate a stochastic process.
We denote the random variables of our process as follow: the gradient estimator $ G_k$, the hessian estimator $\mathcal H_k$, the step size parameter $ \mathcal T_k$, the search direction $\cS_k$, the iterate $X_k$. Their realizations are denoted as  $g_k = G_k(\omega_k), H_k = \mathcal H_k(\omega_k), t_k = \mathcal T_k(\omega_k)$, $s_k = \cS_k(\omega_k)$ and $x_k=X_k(\omega_k)$, respectively, with $\omega_k$ taken from a proper probability space. For brevity we will omit  $\omega_k$ in the following.
We let $\mathcal{E}_{k-1}$ denote all noise history up to iteration $k - 1$ and we include
$\mathcal{E}_{k-1}$ in the algorithmic history for completeness. 
We use $\mathcal{F}_{k-1} = \sigma(G_0,\ldots, G_{k-1}, \mathcal H_0,\ldots ,\mathcal H_{k-1},\mathcal{E}_{k-1})$
to denote the $\sigma$-algebra generated by $G_0, \ldots, G_{k-1}, \mathcal H_0,\ldots ,\mathcal H_{k-1}$ and $\mathcal{E}_{k-1}$, 
up to the beginning  of iteration $k$.

We assume that the random gradient estimators $G_k$ are $(1-\delta_g)$-sufficiently accurate.
\begin{ass}{\rm (gradient estimate)} \label{gacc}The estimator $G_k$ is $(1-\delta_g)$-probabilistically sufficient accurate
in the sense that the indicator variable 
$$
I_k = \1 \{ \|  G_k - \nabla f(X_k) \| \leq {t_k\eta_k} \| G_k \| \}
$$
satisfies the  submartingale condition
\be
 P\left( I_k=1 | \mathcal{F}_{k-1}  \right) \geq 1 - \delta_g, \quad  \delta_g\in(0,1).
\ee
\end{ass}
\vskip 5pt
Iteration $k$ is  \emph{true} iteration if  $I_k=1$, \emph{false} otherwise. Trivially,  if $k$th iteration is true then  the triangle inequality implies
\be \label{truecons}
\| \nabla f(x_k) \| \leq \left( {1+t_k\eta_k} \right) \| g_k \|.
\ee

\vskip 5pt
Finally, at each iteration and for any realization,  the approximation $H_k$ is supposed to be positive definite.

\begin{ass} \label{assH} For all $k\ge 0$, $H_k$ is symmetric positive definite and 
\be \label{convexityH} 
\lambda_1 I \preceq  H_k \preceq \lambda_n I, 
\ee
with $\lambda_1, \, \lambda_2$ as in (\ref{convexity}).
\end{ass}

We remark that assuming $f$ strictly convex, (\ref{convexity}) and (\ref{convexityH}) hold,
with suitable choices of the scalar $\lambda_1, \, \lambda_n$,
as long as the sequence  $\{H_k\}$ is symmetric positive definite and has eigenvalues 
uniformly bounded from below and above.

{Some comments on Assumptions \ref{ass1f}, \ref{ass2f}, \ref{gacc} and \ref{assH} are in order.
They appear in a series of papers on unconstrained optimization where the evaluation of function, gradient  
and Hessian are inexact, either with a controllable or a random noise. Controlling the noise 
in a deterministic way, as in Assumptions \ref{ass1f} and \ref{ass2f} is a realistic request
in applications such as those where the accuracy of $f$-values can be enforced by the magnitude of some discretization parameters or 
$f$ is approximated by using smoothing operators \cite{MWDS,MS,NS}. Probabilistically sufficient accurate gradients,
as in Assumption (\ref{gacc}), 
occur when the gradients are estimated by finite difference and some computation fails to complete, 
in derivative-free optimization or when the gradient are estimated by sample average approximation methods 
 \cite{belgmt,BSV, cartis, STORM1}.
Finally, Assumption \ref{assH} amounts to building a convex random model and is trivially   enforced if the  Hessian
is  estimated by sample average approximation methods and Assumption \ref{convexity} holds.
In literature, unconstrained optimization with inexact function and derivative evaluations 
covers many cases: exact function and gradient evaluations and possibly random Hessian \cite{BGM, noc2, xu2},
exact function and random gradient and Hessian \cite{BSV,bg, cartis},
approximated function and random gradient and Hessian \cite{BellGuriMoriToin22, Berahas}, 
 random function gradient and Hessian \cite{bkm, STORM2, STORM1,PS}. 
In the class of Inexact Newton method with random models  we mention \cite{bkk, bbn, noc1, LR}.

}

\section{{Bounded noise on $f$}}\label{sec3}
{In this section we present and analyze  an Inexact Newton method with line-search where 
the function evaluation is noisy in the sense of Assumption \ref{ass1f}.

At iteration $k$, given $x_k$ and the steplength $t_k$, a  non-monotone Armijo condition  given in \cite{Berahas}
is used. It employs the known upper bound $\varepsilon_f$ introduced in Assumption \ref{errf} and has the form
\be \label{armijo}
\tilde f(x_k + t_k s_k) \leq \tilde f(x_k) + c t_k s_k^T  g_k + 2\varepsilon_f
\ee
$c\in(0,1)$.  If $x_k+t_ks_k$ satisfies (\ref{armijo}) we say that the iteration is successful, we accept the step and increase the step-length $t_k$ for the next iteration. Otherwise the step is rejected and the step-length $t_k$ is reduced. Our procedure belongs to the framework given in Section 4.2 of \cite{Berahas} and it is sketched in Algorithm \ref{algo}.
}

 \algo{}{$k$-th iteration of Algorithm}{
\noindent
Given  $ x_k\in {\mathbb R}^n, \; c \in (0,1)$, $\bar \eta \in (0,1)$ , $\tau \in (0,1)$,  $t_{\max}>0$, $t_k \in(0, t_{\max}]$.
\vskip 2 pt
\begin{itemize}
\item[Step 1.] Choose $\eta_k\in (0,\bar \eta)$.
\item[Step 2.] Apply  CG method with null initial guess to  $ H_k s=-  g_k$ and compute
$s_k$ satisfying (\ref{inexact1}).
\item[Step 3.] If $t_k$ satisfies condition (\ref{armijo}) (successful iteration) then \\
$ x_{k+1} = x_{k} + t_k s_k$,\quad  $t_{k+1}=\max\{t_{\max},\tau^{-1} t_k\}$,\quad $ k = k+1$.
\item[] else $ x_{k+1} = x_{k}$,\quad $t_{k+1}=\tau t_k$,\quad $ k = k+1$.
\end{itemize}
}\label{algo}

{
Similarly to \cite{Berahas, cartis}, Algorithm \ref{algo} generates a stochastic process. 
Given $x_{k}$ and $t_{k}$, the iterate $x_{k+1}$   is fully determined by $g_{k}$, $H_k$ and the noise in
the function value estimation during iteration $k $.
}

{Concerning the well definiteness of the linesearch strategy  (\ref{armijo}), we now prove that if the iteration is true
and $t_k$ is small enough, the linesearch condition is satisfied.}
\vskip 5pt

\begin{lem}\label{lowertk} 
{Suppose that Assumptions \ref{strconv}, \ref{ass1f}, \ref{gacc} and \ref{assH} hold. Suppose that iteration $k$ is true 
and consider any realization of Algorithm \ref{algo}.
Then the iteration is successful whenever $t_k \leq \bar{t}= \frac{2\beta \kappa_1(1-c)}{\kappa_1\lambda_n+2}$.}
\end{lem}

\begin{proof}
 Let $ k $ be an arbitrary iteration.  
Assumptions \ref{strconv} and (\ref{errf}) imply, using the standard arguments for functions with bounded Hessians,
\begin{eqnarray*}
\tilde f(x_k + t_ks_k) - \varepsilon_f \leq f(x_{k}+t_{k} s_{k})&=&
f(x_{k}) +  \int_{0}^{1} [\nabla f(x_{k}+y t_{k} s_{k})]^{T} ( t_{k}s_{k}) dy\\
& = & f(x_{k}) +  \int_{0}^{1} t_k\left([\nabla f(x_{k}+y t_{k} s_{k})]^{T} s_{k} \pm \nabla f(x_k)^Ts_k\right)dy\\
& = & f(x_{k}) +  \int_{0}^{1} t_k[\nabla f(x_{k}+y t_{k} s_{k}) - \nabla f(x_k)]^{T} s_{k} dy + t_k\nabla f(x_k)^Ts_k\\
& \leq & f(x_{k}) +  \int_{0}^{1} t_k\|\nabla f(x_{k}+y t_{k} s_{k}) - \nabla f(x_k)\| \| s_{k}\| dy + t_k\nabla f(x_k)^Ts_k\\
& \leq & f(x_{k}) +  \frac{\lambda_n}{2}  t_{k}^2 \|s_{k}\|^2 + t_k\nabla f(x_k)^Ts_k\\
& \leq & \tilde f(x_{k}) + \varepsilon_f +  \frac{\lambda_n}{2}  t_{k}^2 \|s_{k}\|^2 + t_k\nabla f(x_k)^Ts_k.
\end{eqnarray*}

Since iteration $k$ is true by assumption then  $\|\nabla f(x_k)- g_k\|\le t_k \eta_k \|g_k\|$ holds, and  
using Lemma \ref{Lemma3a} we obtain
\begin{eqnarray}
\tilde f(x_k+t_ks_k) &\leq& \tilde f(x_{k}) + 2\varepsilon_f +  \frac{\lambda_n}{2}  t_{k}^2 \|s_{k}\|^2 + t_k \nabla f(x_k)^Ts_k \pm t_k  g_k^Ts_k \nonumber\\
&=& \tilde f(x_{k}) + 2\varepsilon_f +  \frac{\lambda_n}{2}  t_{k}^2 \|s_{k}\|^2 + t_k [\nabla f(x_k)- g_k]^Ts_k + t_k  g_k^Ts_k \nonumber\\
&\leq& \tilde f(x_{k}) + 2\varepsilon_f +  \frac{\lambda_n}{2}  t_{k}^2 \|s_{k}\|^2 + t_k^2 \frac{\eta_k}{\kappa_1} \|s_k\|^2 + t_k  g_k^Ts_k. \label{disugdim}
\end{eqnarray}

Then, the linesearch condition (\ref{armijo}) is clearly enforced whenever
\begin{eqnarray*}
\tilde f(x_{k}) + 2\varepsilon_f +  \frac{\lambda_n}{2}  t_{k}^2 \|s_{k}\|^2 + t_k^2 \frac{\eta_k}{\kappa_1} \|s_k\|^2 + t_k  g_k^Ts_k 
\leq  \tilde f(x_k) + c t_k s_k^T  g_k + 2\varepsilon_f  
\end{eqnarray*}
which gives
\begin{eqnarray*}
t_k \|s_k\|^2 \left( \frac{\lambda_n}{2} + \frac{\eta_k}{\kappa_1} \right) \leq - (1-c)  g_k^Ts_k .
\end{eqnarray*}
Using  (\ref{hpconstant}) we have $-(1-c) g_k^Ts_k \geq(1-c)  \beta \| s_k \|^2 $. 
Then, since $\eta_k<\bar \eta<1$,  if 
\begin{eqnarray*}
t_k \|s_k\|^2 \left( \frac{\lambda_n}{2} + \frac{1}{\kappa_1} \right) \leq (1-c) \beta \| s_k \|^2,
\end{eqnarray*}
then (\ref{armijo}) holds
and this yields the thesis.
\end{proof}

\subsection{Complexity analysis of the stochastic process}
{ In this section we carry out the convergence analysis of Algorithm \ref{algo}. To this end
we provide a bound on the expected number of
iterations that the algorithm  takes before it achieves a desired level of accuracy in the optimality gap
$f(x_k)-f^*$ with $f^*=f(x^*)$ being the minimum value attained by $f$.
Such a number  of iteration is defined formally below.

\begin{defi}\label{ht}
Let $x^*$ be the global minimizer of $f$ and $f^* = f(x^*)$. Given some 
$\epsilon>0$, $N_\epsilon$ is the number of iterations required until $f(x_k) - f^* \leq \epsilon$ occurs for the first time. 
\end{defi}
The number of iterations $N_\epsilon$ is a random variable and it can be defined as the hitting time for our stochastic process. Indeed it has the property $\sigma(\1 \{N_\epsilon > k\}) \subset \mathcal{F}_{k-1}$.

Following  the notation introduced in Section \ref{sec2}
we let $X_k$, $k\ge 0$, be the random variable with realization $x_k=X_k(\omega_k)$ and  consider
the  following measure of progress towards optimality: 
\be \label{zk}
Z_k= \log \left( \frac{f(X_0)-f^*}{f(X_k)-f^*} \right).
\ee
Further, we let 
\be
Z_\epsilon = \log \left( \frac{f(X_0)-f^*}{\epsilon} \right), \label{ze}
\ee
be  an upper bound for $Z_k$ for any $k<N_\epsilon$. We denote with
$z_k=Z_k(\omega_k)$  a realization of the random quantity $Z_k$.
}

A theoretical framework for analyzing    a generic line search with noise has been developed
in \cite{Berahas}. Under a suitable set of conditions, it provides the expected value for $N_\epsilon$.
We state a result from \cite{Berahas} and will exploit it for our algorithm.

\begin{teo} \label{EN}
Suppose that Assumptions \ref{strconv}, \ref{ass1f}, \ref{gacc}, \ref{assH} hold.
Let $z_k$ a realization of $Z_k$ in (\ref{zk}) and suppose that there exist a constant $\bar t > 0$, a nondecreasing function $h(t) : \RR^+\rightarrow\RR $, which satisfies $h(t) > 0$ for any $t \in(0,t_{max}]$, and a nondecreasing function $r(\varepsilon_f ) : \RR\rightarrow \RR$, which satisfies $r(\varepsilon_f ) \geq 0$ for any $\varepsilon_f \geq 0$, such that for any realization of Algorithm \ref{algo} the following hold for all $k < N_\epsilon$:
\begin{enumerate}
\item[(i)] If iteration $k$ is true and successful, then $z_{k+1} \geq z_k +h(t_k)- r(\varepsilon_f )$.
\item[(ii)] If $t_k \leq \bar t$ and iteration $k$ is true then iteration $k$ is also successful, which implies $t_{k+1} = \tau^{-1}t_k$.
\item[(iii)] $z_{k+1} \geq z_k - r(\varepsilon_f )$ for all successful iterations $k$ and $z_{k+1} \geq z_k$ for all unsuccessful iteration $k$.
\item[(iv)] The ratio $r(\varepsilon_f )/h(\bar t)$ is bounded from above by some $\gamma \in (0, 1)$.
\end{enumerate}
Then under the condition that $\delta_g < \frac{1}{2} - \frac{\sqrt{\gamma}}{2} $, the stopping time $N_\varepsilon$ is bounded in expectation as follows
\be
\E[N_{\epsilon}] \leq \frac{2(1-\delta_g)}{(1-2\delta_g)^2-\gamma}\left[  \frac{2Z_\epsilon}{h(\bar t)} + (1-\gamma)\log_{\tau}\frac{\bar t}{t_0} \right]
\ee
\begin{proof}
See \cite[Assumption 3.3 and Theorem 3.13]{Berahas}.
\end{proof}
\end{teo}
\vskip 5pt
{ We show that our algorithm satisfies the assumptions in Theorem \ref{EN}
if the magnitude of $\epsilon$ fulfills the following condition.

\begin{ass}  \label{rapp} Let $c\in(0,1)$ as in (\ref{armijo}), $\beta$, $\kappa_1$ as in Lemma \ref{Lemma3a}, $\lambda_1$ as in Assumption \ref{strconv}, $\bar t$ as in Lemma \ref{lowertk}, $t_{max}$ as in Algorithm \ref{algo}. 
Assume that $\epsilon$ in Definition \ref{ht} is such that
\be
\epsilon > \frac{4\varepsilon_f}{\left(1- M\right)^{-\gamma}-1}
\ee 
where $M=\frac{c\beta\kappa_1^2\lambda_1 \bar t}{(1+t_{\max})^2}$, for some $\gamma\in(0,1)$.
\end{ass}
{Note that the smaller $M$ is, the larger is $\epsilon$ with respect to $\varepsilon_f$.}

First, we provide a relation between
$z_k$ and $z_{k+1}$ of the form specified in item (i), Theorem \ref{EN}.}

\begin{lem} \label{hrstrconv} 
{Suppose that Assumptions \ref{strconv}, \ref{ass1f}, \ref{gacc}, \ref{assH}, \ref{rapp} hold.}
Consider any realization of Algorithm \ref{algo}. If the $k$-th iterate is true and successful, then 
\be
z_{k+1} \geq z_k - \log \left(1-c \beta\kappa_1^2\lambda_1 \frac{t_k}{(1+t_{\max})^2}\right) - \log\left(1 + \frac{4\varepsilon_f}{\epsilon}\right)
\ee
whenever $k<N_{\epsilon}$.
\begin{proof}
By (\ref{conv1}), $f(x_k)-f^* \leq \frac{1}{2\lambda_1}\| \nabla f(x_k) \|^2$. Using (\ref{truecons}) 
\be
\| g_k \|^2 \geq \left( \frac{1}{1+t_k\eta_k} \right)^2 \| \nabla f(x_k)  \|^2 \geq 2\lambda_1  \left(  \frac{1}{1+t_k\eta_k}\right)^2 ( f(x_k)-f^*) \label{tri}
\ee
Combining condition (\ref{armijo}),  (\ref{errf})  and Lemma \ref{Lemma3a}  it holds
\begin{eqnarray}
 f(x_k) -  f(x_{k+1})  + 2 \varepsilon_f \geq \tilde f(x_k) - \tilde f(x_{k+1}) &\geq& -ct_k g_k^T s_k - 2 \varepsilon_f \label{3*}\\
 &\geq& ct_k \beta\kappa_1^2  \| g _k\|^2 - 2 \varepsilon_f \nonumber
\end{eqnarray}
and thus, using \eqref{tri},
$$
 f(x_k) -  f(x_{k+1}) \pm f^*\geq  2 ct_k  \beta\kappa_1^2\lambda_1   \left(\frac{1}{1+t_k\eta_k}\right)^2 ( f(x_k)-f^*) - 4 \varepsilon_f.
$$
Then it holds
$$
 f(x_{k+1}) -f^* \leq  \left(1-2 c  \beta\kappa_1^2\lambda_1  \frac{t_k}{(1+t_k\eta_k)^2} \right)( f(x_k)-f^*) + 4 \varepsilon_f.
$$
We define $\Delta_k^f = f(x_k) - f^*$. Because of $f(x_k)-f^*>\epsilon$ we have
\begin{eqnarray*}
\Delta_{k+1}^f &\leq&  \left(1-2 c  \beta\kappa_1^2\lambda_1 \frac{t_k}{(1+t_k\eta_k)^2}  + \frac{4\varepsilon_f}{\epsilon}\right) \Delta_k^f\\
&\leq&  \left(1- c  \beta\kappa_1^2\lambda_1 \frac{t_k}{(1+t_k\eta_k)^2}  -\frac{4\varepsilon_f}{\epsilon}c  \beta\kappa_1^2\lambda_1  \frac{t_k}{(1+t_k\eta_k)^2}  + \frac{4\varepsilon_f}{\epsilon}\right) \Delta_k^f\\
&=& \left(1-c \beta\kappa_1^2\lambda_1 \frac{t_k}{(1+t_k\eta_k)^2}  \right)\left( 1 + \frac{4\varepsilon_f}{\epsilon}\right) \Delta_k^f\\
&\leq& \left(1-c  \beta\kappa_1^2\lambda_1\frac{t_k}{(1+t_{\max})^2} \right)\left( 1 + \frac{4\varepsilon_f}{\epsilon}\right) \Delta_k^f\\
\end{eqnarray*}
where the second inequalities holds thanks to Assumption \ref{rapp}, because $4\varepsilon_f <\epsilon$ and the last one holds since $t_k\leq t_{\max}$ and $\eta_k\leq1$.\\
Notice that since $\left( 1+ \frac{4\varepsilon_f}{\epsilon}\right) >~0$, $\Delta_k^f > 0$ and $\Delta_{k+1}^f \geq 0$, it holds
$\left(1-c  \beta\kappa_1^2\lambda_1 \frac{t_k}{(1+t_{\max})^2}\right)\geq~0$. Now taking the inverse and then the $\log$ of both sides, adding $\log \Delta_0^f$,
we have
$$
\log \left( \frac{\Delta_0^f}{\Delta_{k+1}^f} \right) \geq \log \left( \frac{\Delta_0^f}{\Delta_k^f} \right) - \log \left(1-c  \beta\kappa_1^2\lambda_1 \frac{t_k}{(1+t_{\max})^2}\right) - \log \left( 1 + \frac{4\varepsilon_f}{\epsilon}\right)
$$
which completes the proof.

\end{proof}
\end{lem}

\vskip 5pt
{ The next lemma analyze item (ii) of Theorem \ref{EN}}
\begin{lem} \label{falsez}
Suppose that Assumptions \ref{strconv}, \ref{ass1f}, \ref{gacc}, \ref{assH} hold.
Consider any realization of Algorithm \ref{algo}. For every iteration that is false and successful, we have
$$
z_{k+1} \geq z_k - \log \left( 1+\frac{4\varepsilon_f}{\epsilon} \right)
$$
Moreover $z_{k+1}=z_k$ for any unsuccessful iteration.
\begin{proof}
For every false and successful iteration, using (\ref{errf}),  (\ref{armijo}) and (\ref{hpconstant}) we have 
\begin{eqnarray*}
f(x_{k+1}) &\leq& f(x_k) + ct_k s_k^\top g_k + 4\varepsilon_f \\
&\leq & f(x_k) + 4\varepsilon_f 
\end{eqnarray*}
thus, because of $f(x_k)-f^*>\epsilon$,
\begin{eqnarray*}
f(x_{k+1}) -f^* &\leq&  f(x_k) - f^* + 4\varepsilon_f \\
& \leq & \left( 1+\frac{4\varepsilon_f}{\epsilon} \right)(f(x_k) - f^*).
\end{eqnarray*}
So it holds $\Delta_{k+1}^f \leq \left( 1+\frac{4\varepsilon_f}{\epsilon}  \right)\Delta_k^f$.
Now taking the inverse and then the log of both sides, adding $\log \Delta_0^f$ we have
$$
\log \left( \frac{\Delta_0^f}{\Delta_{k+1}^f} \right) \geq \log \left( \frac{\Delta_0^f}{\Delta_k^f} \right)  - \log \left( 1 + \frac{4\varepsilon_f}{\epsilon}\right)
$$
which completes the first part of the proof. Finally for any unsuccessful iteration $z_{k+1}=z_k$ follows by Step 3 of Algorithm \ref{algo} that provides $x_{k+1}=x_k$ and hence $f(x_{k+1})=f(x_k)$ whenever iteration $k$ is not successful. 
\end{proof}
\end{lem}
\vskip 5pt
{ We can now summarize our results. First, note that $\left(1-c  \beta\kappa_1^2\lambda_1 \frac{t_k}{(1+t_{\max})^2}\right) \geq 0$ for all $t_k\in[0,t_{max}]$, due to (\ref{betak1}) and (\ref{k11}).
Second, let 
\be
h(t) = - \log \left(1-c  \beta\kappa_1^2\lambda_1\frac{t}{(1+t_{\max})^2}\right) \quad \text{and} \quad r(\varepsilon_f) =  \log \left( 1 + \frac{4\varepsilon_f}{\epsilon}\right).
\ee
It is easy to see that $h(t)$ is monotone and non increasing if $t\in [0,t_{\max}]$.

Combining Lemma \ref{lowertk}, \ref{hrstrconv} and \ref{falsez}, 
we have that for any realization of Algorithm \ref{algo} and $k< N_\epsilon$ with $\epsilon$ as  in Assumption \ref{rapp}:}
\begin{enumerate}
\item[(i)] (Lemma \ref{hrstrconv}) If iteration $k$ is true and successful, then $z_{k+1} \geq z_k +h(t_k)- r(\varepsilon_f )$.
\item[(ii)] (Lemma \ref{lowertk}) If $t_k \leq \bar t$ and iteration $k$ is true then iteration $k$ is also successful, which implies $t_{k+1} = \tau^{-1}t_k$.
\item[(iii)] (Lemma \ref{falsez}) $z_{k+1} \geq z_k - r(\varepsilon_f )$ for all successful iterations $k$ and $z_{k+1} = z_k$ for all unsuccessful iteration $k$.
\item[(iv)] (Assumption \ref{rapp}) The ratio $r(\varepsilon_f )/h(\bar t)$ is bounded from above by some $\gamma \in~(0, 1)$.
\end{enumerate}

Hence, we can use Theorem \ref{EN}  and get the following boun on $\E[N_\epsilon]$,
$$
\E[N_\epsilon] \leq \frac{2(1-\delta_g)}{(1-2\delta_g)^2-\gamma}\left[  {2 \log_{1/(1-M)} \left( \frac{f(x_0)-f^*}{\epsilon} \right)} + (1-\gamma)\log_{\tau}\frac{\bar t}{t_0} \right]
$$
with $M$ given in Assumption \ref{rapp}.
{This result is valid under Assumption \ref{rapp},
namely for sufficiently large values of $\epsilon$. The fact that $\epsilon$ cannot be arbitrarily small
is consistent with the presence of
noise $\varepsilon_f$ in $f$-evaluations. Trivially, if $\varepsilon_f=0$ then the optimality gap $f(x_k)-f^*$ can be made arbitrarily small.}

\section{{Decreasing noise on $f$}}\label{sec4}
{
In this section we present an Inexact Newton algorithm suitable to the case where 
$f$-evaluations can be performed with adaptive accuracy.
We use the linesearch condition 
\be \label{armijo2}
\tilde f(x_k + t_k s_k) \leq \tilde f(x_k) + c t_k s_k^T  g_k ,
\ee
with $\tilde f(x_k)$ and $\tilde f(x_k + t_k s_k)$ satisfing Assumption \ref{ass2f}
with $\theta < \frac{c}{2}$.
In fact, (\ref{armijo2}) has the form of the classical Armijo condition but the true $f$ is replaced by the approximation~$\tilde f$.

The resulting algorithm is given below.}

 \algo{}{$k$-th iteration of Algorithm}{
\noindent
Given  $ x_k\in {\mathbb R}^n, \; c \in (0,1)$, $\bar \eta \in (0,1)$ , $\tau \in (0,1)$, 
$\theta>0$,  $t_{\max}>0$, $t_k \in(0, t_{\max}]$.
\vskip 2 pt
\begin{itemize}
\item[Step 1.] Choose $\eta_k\in (0,\bar \eta)$.
\item[Step 2.] Apply  CG method initialized by the null vector to  $ H_k s_k=-  g_k$ and compute
$s_k$ satisfying (\ref{inexact1}).
\item[Step 3.] Compute $\tilde f(x_k)$ and $\tilde f(x_k + t_k s_k)$ satisfying  (\ref{errfk}).
\item[Step 4.] If $t_k$ satisfies condition (\ref{armijo2}) (successful step) then \\
$ x_{k+1} = x_{k} + t_k s_k$,\quad  $t_{k+1}=\max\{t_{\max},\tau^{-1} t_k\}$,\quad $ k = k+1$.
\item[] else $ x_{k+1} = x_{k}$,\quad $t_{k+1}=\tau t_k$,\quad $ k = k+1$.
\end{itemize}
}\label{algo2}

The following Lemma shows that 
a successful iteration is guaranteed whenever it is true and  $t_k$ is sufficiently small.
\vskip 5pt\noindent
\begin{lem} \label{lowertk2} 
{Suppose that Assumptions \ref{strconv}, \ref{ass2f} with $\theta<\frac{c}{ 2}$, \ref{gacc} and \ref{assH} hold.} Suppose that iteration $k$ is true 
and consider any realization of Algorithm \ref{algo2}.
Then the iteration is successful whenever  $t_k\leq \bar t = \frac{2\kappa_1\beta}{\kappa_1\lambda_n+2}(1-c-2\theta)$.
\begin{proof}
Using the same arguments as in Lemma \ref{lowertk}, using (\ref{errfk}) and (\ref{armijo2}), rather than (\ref{errf}) and (\ref{armijo}) we obtain
\begin{eqnarray*}
\tilde f(x_k + t_ks_k) 
&\leq& \tilde f(x_{k}) - 2\theta t_ks_k^Tg_k  +  \frac{\lambda_n}{2}  t_{k}^2 \|s_{k}\|^2 + t_k^2 \frac{\eta_k}{\kappa_1} \|s_k\|^2 + t_k  g_k^Ts_k.
\end{eqnarray*}
The linesearch condition (\ref{armijo2}) is clearly enforced whenever
$$
\tilde f(x_{k}) - 2\theta t_ks_k^Tg_k  +  \frac{\lambda_n}{2}  t_{k}^2 \|s_{k}\|^2 + t_k^2 \frac{\eta_k}{\kappa_1} \|s_k\|^2 + t_k  g_k^Ts_k \leq  \tilde f(x_k) + c t_k s_k^T  g_k
$$
which gives
\begin{eqnarray*}
t_k \|s_k\|^2 \left( \frac{\lambda_n}{2} + \frac{\eta_k}{\kappa_1} \right) \leq - (1-c-2\theta)  g_k^Ts_k 
\end{eqnarray*}

Using  (\ref{hpconstant}) we have $-(1-c-2\theta) g_k^Ts_k \geq(1-c-2\theta)  \beta \| s_k \|^2 $. 
Then, since $\eta_k<\bar \eta<1$, (\ref{armijo2}) holds if 
\begin{eqnarray*}
t_k \|s_k\|^2 \left( \frac{\lambda_n}{2} + \frac{1}{\kappa_1} \right) \leq (1-c-2\theta) \beta \| s_k \|^2 
\end{eqnarray*}

and this yields the thesis.
\end{proof}
\end{lem}

\subsection{Complexity analysis of the stochastic process}
The behaviour of the method is studied analyzing the hitting time $N_\epsilon$ in 
Definition \ref{ht}.  In particular, we first show the following two results on
the realization $z_k$ of the variable $Z_k$ in (\ref{zk}).
\begin{lem} \label{hrstrconv2} 
{Suppose that Assumptions \ref{strconv}, \ref{ass2f} with $\theta<\frac{c}{ 2}$, \ref{gacc} and \ref{assH} hold.   
If the $k$-th iterate of Algorithm \ref{algo2} is true and successful, for any realization of the Algorithm \ref{algo2} we have}
\be
z_{k+1} \geq z_k - \log \left(1-  2 \beta \kappa_1^2 \lambda_1  \left( c-2\theta \right)  \frac{t_k}{(1+t_{max})^2} \right),
\ee
whenever $k<N_{\epsilon}$.
\begin{proof}
Using the same arguments as in Lemma \ref{hrstrconv}, using (\ref{errfk}) and (\ref{armijo2}), rather than (\ref{errf}) and (\ref{armijo}) we obtain
\begin{eqnarray*}
 f(x_{k}) -f(x_{k+1})  &\geq& -(c-2\theta)t_k g_k^Ts_k\\
 &\geq& t_k \beta 	\kappa_1^2\left( c-2\theta \right) \| g_k \|^2,
 \end{eqnarray*}
 where the second inequality comes from  (\ref{hpconstant}). 
 Thus
 $$
 f(x_{k+1})-f^* \leq  f(x_k) -f^* -t_k\beta\kappa_1^2(c-2\theta) \| g_k \|^2
 $$
 and using (\ref{tri}) we get
 $$
f(x_{k+1})-f^* \leq  \left(1-  2 \beta \kappa_1^2 \lambda_1  \left( c-2\theta \right)  \frac{t_k}{(1+t_{max})^2}  \right)( f(x_k)-f^*).
$$
Now proceeding as in Lemma \ref{hrstrconv} we have the thesis. 
\end{proof}
\end{lem}

\begin{lem} \label{falsez2}
{Suppose that Assumptions \ref{strconv}, \ref{ass2f} with $\theta<\frac{c}{ 2}$, \ref{gacc} and \ref{assH} hold.   
For any realization of Algorithm \ref{algo2}  we have
$$
z_{k+1} \geq z_k,
$$
if  the iteration $k$ is false and successful, 
$$z_{k+1}=z_k,$$
if the iteration $k$ is  unsuccessful.}
\begin{proof}
For every false and successful iteration, using (\ref{errfk}) and (\ref{armijo2}),we have 
\begin{eqnarray*}
f(x_{k+1}) &\leq& f(x_k) + ct_k s_k^\top g_k - 2\theta t_kg_k^Ts_k \\
&=& f(x_k) + (c-2\theta) t_k s_k^Tg_k\\
&\leq & f(x_k),
\end{eqnarray*}
and  in case of unsuccessful iteration Step 4 of the algorithm provides $x_{k+1}=x_k$.

\end{proof}
\end{lem}

Now we can state the main result on the expected value for the hitting time.

\begin{teo} \label{EN2}
 Suppose that Assumptions \ref{strconv}, \ref{ass2f} with $\theta<\frac{c}{2}$, \ref{gacc} and \ref{assH} hold
 and let $\bar t$ given in Lemma \ref{lowertk2}.
Then under the condition that $\delta_g < \frac{1}{2} $, the stopping time $N_\varepsilon$ is bounded in expectation as follows
$$
\E[N_\epsilon] \leq \frac{2(1-\delta_g)}{(1-2\delta_g)^2}\left[  {2 \log_{1/(1-M)} \left( \frac{f(x_0)-f^*}{\epsilon} \right)} + \log_{\tau}\frac{\bar t}{t_0} \right]
$$
with $M=\frac{2(c-2\theta)\beta\kappa_1^2\lambda_1 \bar t}{(1+t_{\max})^2}$, {$\bar t$ as in Lemma (\ref{lowertk2}).}
\begin{proof}
Let  
\be
h(t) = - \log \left(1-  2 \beta \kappa_1^2 \lambda_1  \left( c-2\theta \right)  \frac{t}{(1+t_{max})^2} \right),
\ee
and note that $h(t)$ is non decreasing for $t\in [0,t_{\max}]$ and that $h(t)> 0$  for $t\in [0,t_{\max}]$.
For any realization $z_k$ of $Z_k$ in (\ref{zk}) of Algorithm \ref{algo2} 
the following hold for all $k < N_\epsilon$:
\begin{enumerate}
\item[(i)] If iteration $k$ is true and successful, then $z_{k+1} \geq z_k +h(t_k)$ by Lemma \ref{hrstrconv2}.
\item[(ii)] If $t_k \leq \bar t$ and iteration $k$ is true then iteration $k$ is also successful, which implies $t_{k+1} = \tau^{-1}t_k$ by Lemma  \ref{lowertk2}.
\item[(iii)] $z_{k+1} \geq z_k $ for all successful iterations $k$ ($z_{k+1} = z_k$ for all unsuccessful iteration $k$), by Lemma \ref{falsez2}.
\end{enumerate}
Moreover,  our stochastic process $\{{\mathcal{T}_k}, Z_k\}$ obeys the expressions below.
By Lemma \ref{lowertk2} and the definition of Algorithm \ref{algo2} the update of the random variable $\mathcal T_k $ such that $t_k = \mathcal T_k(\omega_k)$  is
$$
\mathcal T_{k+1}=  \left\{
\begin{array}{ll}
\tau^{-1} \mathcal T_k & \mbox{ if } I_k=1 , \, \mathcal T_k \le \bar t \  \mbox{ (i.e., successful)}\\
\tau^{-1} \mathcal T_k & \mbox{ if the iteration is successful}, \,  I_k=0 , \, \mathcal T_k \le \bar t \\
\tau\, \mathcal T_k & \mbox{ if the iteration is unsuccessful}, \,  I_k=0 , \, \mathcal T_k \le \bar t \\
\tau^{-1} \mathcal T_k & \mbox{ if the iteration is successful},\mathcal T_k > \bar t \\
\tau\, \mathcal T_k & \mbox{ if the iteration is unsuccessful},   \mathcal T_k > \bar t \\
\end{array} 
\right.
$$
By  Lemma \ref{lowertk2}  Lemma \ref{hrstrconv2} and  Lemma \ref{falsez2} the  random variable $Z_k $ obeys the expression
$$
Z_{k+1}\ge   \left\{
\begin{array}{ll}
Z_k +h(\mathcal{T}_k)& \mbox{ if the iteration is successful and }  \,  I_k=1 \\
Z_k & \mbox{ if the iteration is successful and }  I_k=0  \\
Z_k& \mbox{ if the iteration is unsuccessful}\\
\end{array} 
\right.
$$
Then Lemma 2.2--Lemma 2.7 and Theorem 2.1 in \cite{cartis} hold which gives the thesis.
\end{proof}
\end{teo}

\subsection{Local convergence}
We conclude our study analyzing the local behaviour of the Newton-CG method employing gradient estimates $(1-\delta_g)$-probabilistically sufficiently accurate, i.e. satisfying 
Assumption \ref{gacc} and Hessian estimates satisfying the following assumption.
\vskip 5pt
\begin{ass} \label{ascl2}
The Hessian of the objective function $ f$ is Lipschitz-continuous with constant $L_H>0$, 
\be \label{hlip2}
\| \nabla ^2 f(x)-\nabla ^2 f(y) \| \leq L_H \| x-y \|, \quad \forall x, y \in\RR^n
\ee
and there exists a constant $C>0$ such that the Hessian estimator is $(1-\delta_H)$-probabilistically sufficiently accurate as follows
\be \label{vhbound2}
P( J_k=1 | \mathcal{F}_{k-1} ) \geq 1-\delta_H, \quad  \delta_H\in(0,1).
\ee
where $J_k=\1\{ \| \mathcal H_k - \nabla ^2 f(X_k)\|  \leq C\eta_k\}$.
\end{ass}
We let $t_{\max}=1$, so that the maximum step-size gives the full CG step $s_k$. 

{The following lemma shows that if the full CG step $s_k$ is accepted then the error 
linearly decreases with a certain probability. Further, the same occurrence over $\ell$ successive 
iterations is analyzed.}

\begin{lem} \label{convloc2}
Suppose that Assumptions \ref{strconv}, \ref{ass2f} with $\theta<\frac{c}{2}$, \ref{gacc} and \ref{assH} hold. Let $x_{\bar k}$ be a realization of Algorithm \ref{algo2} with $t_{\bar k}=1$.
Assume that  the iteration is successful and  $\|x_{\bar k}-x^*\|$ and $\eta_{\bar k}$ are sufficiently small so that 
$\frac{1}{\lambda_1}\left[ \frac{L_H}{2}  \|x_{\bar k}-x^*\| + C \eta_{\bar k}  +  \frac{2\lambda_n\eta_{\bar k}}{1-\bar \eta} \right] <~\tilde C<~1$. Then, at least with probability $p=(1-\delta_g)(1-\delta_H)$, it holds
$$
\| x_{\bar k+1}-x^* \| < \tilde C \| x_{\bar k}-x^* \|.
$$
If $\{ \eta_k \}$ is a non-increasing sequence and  the iterations $\bar k,\ldots,\bar k+\ell-1$ are successful with $t_k=1$  
 for $k=\bar k,\ldots,\bar k+\ell-1$,
then  it holds $\| x_{ k+1} -x^* \|< \| x_{ k} -x^* \|$ for $k=\bar k,\ldots,\bar k+\ell-1$, at least with probability $p^l$.
\begin{proof}
\begin{eqnarray*}
\| x_{k+1}-x^* \| &=& \| x_k + s_k-x^* \| \\
&=& \| x_k - H_k^{-1}g_k +  H_k^{-1}r_k -x^* \|\\
&=& \|  H_k^{-1} [\nabla^2 f(x_k)(x_k-x^*) - \nabla^2 f(x_k)(x_k-x^*) +  H_k(x_k-x^*) - g_k  \pm \nabla f(x_k)+ r_k  ]\|\\
&\leq& \|  H_k^{-1}\| ( \|\nabla^2 f(x_k)(x_k-x^*) -\nabla f(x_k)\| + \\
&&\| (\nabla^2 f(x_k) - H_k) (x_k-x^*)\| + \| g_k -\nabla f(x_k)\| + \|r_k\| ) \\
&\leq&  \frac{1}{\lambda_1} \Big( \|\nabla^2 f(x_k)(x_k-x^*) -\nabla f(x_k)\| + \\
&& \| (\nabla^2 f(x_k) - H_k) (x_k-x^*)\| + \| g_k -\nabla f(x_k)\| + \|r_k\| \Big).
\end{eqnarray*}

Thanks to (\ref{hlip2}) it holds
$$
 \|\nabla^2 f(x_k)(x_k-x^*) -\nabla f(x_k)\| \leq  \frac{L_H}{2}  \|x_k-x^*\|^2.
$$
Let us assume that both the events $I_k $ and  $J_k$ are true. Then, $\|g_k-\nabla f(x_k)\|\le \eta_k \|g_k\|$,
\begin{eqnarray*}
\| (\nabla^2 f(x_k) - H_k) (x_k-x^*)\|  &\leq& \| \nabla^2 f(x_k) - H_k\| \|x_k-x^*\| \\
&\leq& C \eta_k \|x_k-x^*\|,
\end{eqnarray*}
and by (\ref{inexact1})
$$
\| g_k-\nabla f(x_k)\| + \| r_k \| \leq 2\eta_k \| g_k \|.
$$
Moreover,
$$
\| g_k \| \leq \| g_k-\nabla f(x_k) \| + \| \nabla f(x_k) \| \leq \eta_k \| g_k \| + \| \nabla f(x_k) \|
$$
i.e. $\| g_k \| \leq \frac{1}{1-\eta_k}\| \nabla f(x_k) \|$.
Then combining with (\ref{Lemma2}) we have
$$
\| g_k \| \leq \frac{1}{1-\eta_k}\| \nabla f(x_k) \| \leq \lambda_n\frac{1}{1-\eta_k}\| x_k -x^* \| \leq \lambda_n\frac{1}{1-\bar \eta}\| x_k -x^* \|.
$$
Therefore
$$
\| g_k-\nabla f(x_k)\| + \| r_k \| \leq \frac{2\eta_k\lambda_n}{1-\bar \eta}\| x_k -x^* \|.
$$
Then, since $P(I_k\cap J_k)\geq p$ it follows
\begin{eqnarray*}
\| x_{k+1}-x^* \| &\leq& \frac{1}{\lambda_1} \left[ \frac{L_H}{2}  \|x_k-x^*\|^2 + C \eta_k \|x_k-x^*\| +  \frac{2\lambda_n\eta_k}{1-\bar \eta} \|x_k-x^*\| \right] \\
\end{eqnarray*}
at least with probability $p$. \\
Therefore, since at iteration $\bar k$, $\frac{1}{\lambda_1}\left[ \frac{L_H}{2}  \|x_{\bar k}-x^*\| + C \eta_{\bar k}  +  \frac{2\lambda_n\eta_{\bar k}}{1-\bar \eta} \right] <~\tilde C<~1$ by assumption, it follows $\| x_{\bar k+1}-x^* \| < \tilde C \| x_{\bar k}-x^* \|$. \\
{ At iteration $\bar k+1$, $t_{\bar k+1}=1$ and  the iteration is successful  by hypothesis. Then, we can repeat the previous arguments and the thesis follows.}
\end{proof}
\end{lem}


\section{Finite sum case}\label{sec5}
In this section we consider the finite-sum minimization problem that arises in machine learning and data analysis:
\be \label{finitesum}
\min_{x\in\RR^n} f(x) = \frac{1}{N} \sum_{i=1}^N f_i(x).
\ee
The objective function $f$ is the mean of $N$ component functions $f_i \colon \RR^n \rightarrow \RR$ and for large values of $N$, the exact evaluation of the function and derivatives might be computationally expensive. We suppose that each $f_i$ is strongly convex.

Following \cite{BSV,bg,cartis}  $f$ is evaluated exactly  
while the approximations $g_k$ and $H_k$ to the gradient and the Hessian respectively satisfy accuracy requirements in probability.

The evaluations of $g_k$ and $H_k$ can be made using subsampling, that means picking randomly and uniformly chosen subsets of indexes $\cN_{g,k}$ and $\cN_{H,k}$ from $\cN=\{ 1,\ldots,N \}$ and 
define 
\be \label{gkhksub}
g_k =\frac{1}{| \cN_{g,k} |} \sum_{i\in\cN_{g,k}} \nabla f_i (x_k), \quad \text{ and } \quad H_k =\frac{1}{| \cN_{H,k} |} \sum_{i\in\cN_{H,k}} \nabla^2 f_i (x_k).
\ee 
If we want $g_k$ and $H_k$ to be probabilistically sufficiently accurate as in Definition \ref{gacc} and in Assumption \ref{ascl2} respectively, 
the sample sizes $|\cN_{g,k} |$ and $|\cN_{H,k} |$ can be determined by using the operator-Bernstein inequality introduced in \cite{tropp}. As shown in \cite{belgmt}, 
$g_k$ and $H_k$  are $(1-\delta_g)$ and $(1-\delta_H)$ -probabilistically sufficiently accurate if  
\begin{eqnarray}
|\cN_{g,k} | &\geq& \min \left\{ N, \frac{4 \kappa_{f,g}(x_k)}{\gamma_{g,k}}  \left(\frac{\kappa_{f,g}(x_k)}{\gamma_{g,k}} +\frac{1}{3} \right) \log\left( \frac{n+1}{\delta_g} \right) \right\},  \label{ngk}\\
|\cN_{H,k} | &\geq& \min \left\{ N, \frac{4 \kappa_{f,H}(x_k)}{ C\eta_k }  \left(\frac{\kappa_{f,H}(x_k)}{ C\eta_k } +\frac{1}{3} \right) \log\left( \frac{2n}{\delta_H} \right) \right\}, \label{nhk}
\end{eqnarray}
where $\gamma_{g,k}$ is an approximation of the required gradient accuracy, namely $\gamma_{g,k}\approx~ t_k\eta_k \| G_k \|$
and under the assumption that, for any $x\in\RR^n$, there exist non-negative upper bounds $\kappa_{f,g}$ and $\kappa_{f,H}$ such that
\begin{eqnarray*}
\max_{i\in\{ 1,\ldots,N \}}\| \nabla f_i(x) \| &\leq \kappa_{f,g}(x), \\
\max_{i\in\{ 1,\ldots,N \}}\| \nabla^2 f_i(x) \| &\leq \kappa_{f,H}(x).
\end{eqnarray*}

A practical version of the procedure is shown in Algorithm \ref{algo3}. {
Gradient approximation requires a loop since the accuracy requirement is implicit;
such a strategy is Step 2 of the following algorithm.}
\algo{}{$k$-th iteration of Algorithm}{
\noindent
Given  $ x_k\in {\mathbb R}^n, \; c,\bar \eta, \tau, \kappa_\gamma \in (0,1)$, $\gamma_0>0$,  $t_{\max}>0$, $t_k \in(0, t_{\max}]$.
\vskip 2 pt
\begin{itemize}
\item[Step 1.] Choose $\eta_k\in (0,\bar \eta)$.
\item[Step 2.] Gradient approximation. Set $i=0$ and initialize $\gamma_k^{(i)}=\gamma_0$.
\begin{itemize}
\item[2.1] compute $g_k$ with (\ref{gkhksub}), (\ref{ngk}), $\gamma_{g,k}=\gamma_k^{(i)}$;
\item[2.2] if $\gamma_k^{(i)} \leq t_k\eta_k \| g_k \|$, go to Step 3;
\item[] else, set $\gamma_k^{(i+1)} = \kappa_\gamma \gamma_k^{(i)}$, $i=i+1$, go to Step 2.1.
\end{itemize}
\item[Step 3.] Compute $H_k$ with (\ref{gkhksub}) and (\ref{nhk}).
\item[Step 4.] Apply  CG method initialized by the null vector to  $ H_k s_k=-  g_k$ and compute
$s_k$ satisfying (\ref{inexact1}).
\item[Step 5.] If $t_k$ satisfies 
\be 
 f(x_k + t_k s_k) \leq  f(x_k) + c t_k s_k^T  g_k .
\ee
then $ x_{k+1} = x_{k} + t_k s_k$,\quad  $t_{k+1}=\max\{t_{\max},\tau^{-1} t_k\}$,\quad $ k = k+1$ 
\\(successful iteration)
\item[] else $ x_{k+1} = x_{k}$,\quad $t_{k+1}=\tau t_k$,\quad $ k = k+1$ 
(unsuccessful iteration).
\end{itemize}
}\label{algo3}

Inexact Newton methods for the finite-sum minimization problems are investigated also in 
  \cite{bkk,bbn,mit1}. In \cite{bkk}  it is analyzed a  linesearch Newton-CG method  where    the objective function and the gradient are approximated by subsampling with increasing  samplesizes determined by  a prefixed rule.
Random estimates of  the Hessian with adaptive accuracy requirements as in Assumption \ref{ascl2} are employed and local convergence results  in the mean square are given.
In \cite{bbn} the local convergence  of Inexact Newton method is studied assuming to use prefixed choice of the sample size used  to estimate by subsampling both  the gradient and the Hessian. 
The paper \cite{mit1}  studies the global as well as local convergence behavior of  linesearch Inexact Newton algorithms, where the objective function is exact and the Hessian and/or gradient are sub-sampled. An high probability analysis of the local and local convergence of the method is given, whereas we prove complexity results is expectation with noise in the objective function.  Moreover,  the estimators $g_k$ and $H_k$   are supposed to be $(1-\delta_g)$ and $(1-\delta_H)$ -probabilistic sufficiently accurate as in  our approach but with different   accuracy requirements.  Predetermined and increasing accuracy requirements are used in  \cite{mit1} rather than the adaptive  accuracy requirements  in  \ref{gacc} and in  Assumption \ref{ascl2}.

{\footnotesize
\section*{\footnotesize Acknowledgment}

INdAM-GNCS partially supported the first and third authors under Progetti di
Ricerca 2021. }


\begin{thebibliography}{99}
\bibitem{BSV} A.S Bandeira, K. Scheinberg, L.N. Vicente, {\em Convergence of trust-region methods based on probabilistic models}, SIAM Journal on Optimization, {24}(3), 1238--1264, 2014.
%
\bibitem{bg}  Bellavia, S., Gurioli G.,  
{\em Stochastic analysis of an adaptive cubic regularization method under inexact gradient evaluations and dynamic Hessian accuracy}, Optimization,  71, 227-261, 2022.

\bibitem{BGM} S. Bellavia, G. Gurioli, B. Morini,  {Adaptive cubic regularization methods with dynamic inexact Hessian information and applications to finite-sum minimization}, {\em IMA Journal of Numerical Analysis}, 41(1), 764-799, 2021.
%
\bibitem{BellGuriMoriToin22}
S. Bellavia, G. Gurioli, B. Morini, and {Ph.}L. Toint.
\newblock Adaptive regularization for nonconvex optimization using inexact function values and randomly perturbed derivatives.
\newblock {\em Journal of Complexity},  68, Article number 101591, 2022.


%
\bibitem{belgmt} S. Bellavia, G. Gurioli, B. Morini, Ph.L. Toint, Adaptive Regularization Algorithms with Inexact Evaluations for Nonconvex Optimization. {\em SIAM Journal on Optimization 29(4), 2281-2915, 2019.}
%
\bibitem{bkk} S. Bellavia, N. Krejic, N. Krklec Jerinkic, Subsampled Inexact Newton methods for minimizing large sums of convex functions,
{\em IMA Journal of Numerical Analysis}, 40, 2309-2341, 2020.
%
\bibitem{bkm}  S. Bellavia, N. Kreji\'c, B. Morini, {Inexact restoration with subsampled trust-region
methods for finite-sum minimization}, {\em Computational Optimization and Applications},  {76}, 701-736, 2020. 
%
\bibitem {bbn_2017}  {A.S.Berahas, R. Bollapragada, J. Nocedal}, An Investigation of Newton-Sketch and Subsampled Newton Methods, {\em Optimization Methods and Software}, 
35,  661-680, 2020. 
%
\bibitem {Berahas}  {A.S. Berahas, L.  Cao, K. Scheinberg}, Global convergence rate analysis of a generic line search algorithm with noise, {\em  SIAM Journal on Optimization,} 2019.
%
%
\bibitem{STORM2} J. Blanchet, C. Cartis, M. Menickelly, K. Scheinberg, {Convergence Rate Analysis of a Stochastic Trust Region Method via Submartingales}, {\em INFORMS Journal on Optimization}, {1},   92-119, 2019. 

%
\bibitem{bbn}  Bollapragada, R., Byrd, R., Nocedal, J., Exact and Inexact Subsampled Newton Methods for
Optimization,  {\em IMA Journal Numerical Analysis, 2018.}
%
%
\bibitem{noc2} {R.H. Byrd, G.M. Chin, J. Nocedal, Y. Wu} Sample size selection in optimization methods for machine learning,
 {\em Mathematical Programming}, 134(1), 127-155, 2012.
%
\bibitem{noc1} {R.H. Byrd , G.M. Chin, W. Neveitt, J. Nocedal}, On the Use of Stochastic Hessian Information in Optimization Methods for Machine Learning,
 {\em SIAM Journal on Optimization, 21(3),  977-995, 2011}.
%
%
\bibitem{C} R.G. Carter, On the global convergence of trust-region algorithms using inexact gradient information, {SIMA Journal of Numerical Analysis}, 28, 251-265, 1991.
%
\bibitem {cartis}  {C. Cartis, K. Scheinberg}, Global convergence rate analysis of unconstrained optimization methods based on probabilistic model, {\em  Mathematical Programming}, 169(2), pp. 337--375, 2017.
%
\bibitem{STORM1} R. Chen, M. Menickelly, K. Scheinberg,
{Stochastic optimization using a trust-region method and random models}, {\em Mathematical Programming}, {169}(2), 447-487, 2018.
%
%
 %
%
 \bibitem{DES}  R.S. Dembo, S.C. Eisenstat,  T. Steinhaug, Inexact Newton method, {\em SIAM Journal on  Numerical Analysis 19(2),  400-409, 1982.}
%
%
\bibitem{dskkv} D. di Serafino, N. Krejić, N. Krklec Jerinkić, M. Viola, LSOS: Line-search Second-Order Stochastic optimization methods for nonconvex finite sums, {\em ArXiv:2007.15966v2}, 2021.
\bibitem{EW} S.C. Eisenstat, H.F. Walker,  Choosing the Forcing Terms in an Inexact Newton Method, {\em SIAM Journal on  Scientific  Computing, 17(1), 16-32, 1996.}
%
\bibitem{fg} K. Fountoulakis, J. Gondzio, A second order method for strongly convex $ \ell_1 - $ regularization problems, {\em Mathematical Programming}, 156, pp. 189-219, 2016.

\bibitem{frz} G. Franchini, V. Ruggiero, L. Zanni, Ritz-like values in steplength selections for stochastic gradient  methods, {\em Soft Computing}, 24, 23  17573-17588, 2020.
%
%
%
%
\bibitem{HS} M.R. Hestenes, E. Stiefel, {Methods of conjugate gradients for solving linear systems}, 
{\em Journal of Research of the National Bureau of Standards}, 49(6), pp. 409-436, 1952.
%
%
%
%
%
%
%
%
\bibitem{LR} Y. Liu, F. Roosta, Convergence of Newton-MR under inexact hessian information, {\em
 SIAM Journal on Optimization}, 31(1), 59-90, 2021.
%
\bibitem{MWDS} A. Maggiar, A. Wachter, I.S. Dolinskaya, J. Staum, {A derivative-free trust-region algorithm for the optimization of functions smoothed via gaussian convolution using adaptive multiple importance sampling}, {\em SIAM Journal on
Optimization}, 28, pp. 1478-1507, 2018.
%
\bibitem{MS} J.J. More, S.M. Wild, {Estimating computational noise}, {\em SIAM Journal
on Scientific Computing}, 33, pp. 1292-1314, 2011.
%
%
%

%
%
\bibitem{nesterov} Y. Nesterov, Introductory Lectures on Convex Optimization: A Basic Course, Springer Science and Media, Vol. 87, 2013.
%
\bibitem{NS} Y. Nesterov, V. Spokoiny, Random gradient-free minimization of convex
functions, {\em Foundations of Computational Mathematics}, 17, pp. 527-566, 2017.
%
\bibitem{NW}  Nocedal, J.,  Wright, S. J., 
{Numerical Optimization},
{Springer Series in Operations Research},
{Springer},
1999.
%
\bibitem{PS}   C. Paquette, K. Scheinberg, {A Stochastic Line Search Method with Expected Complexity Analysis},
{\em SIAM Journal of Optimization}, 30 , pp. 349-376, 2020.
%
%
%
%
\bibitem{mit1} F. Roosta-Khorasani, M.W. Mahoney,
Sub-Sampled Newton Methods, {\em Mathematical  Programming},  174,  293--326, 2019.
%
%
%
\bibitem{tropp} J.A. Tropp,  An Introduction to Matrix Concentration Inequalities. {\em Foundations and Trends in Machine Learning 8(1-2), 1--230, 2015.}
%
%
\bibitem{xu2} P. Xu, F. Roosta-Khorasani, M.W. Mahoney, Newton-Type Methods for Non-Convex Optimization Under Inexact Hessian Information,  {\em Mathematical Programming}, 184,  35-70, 2020.
%

\end{thebibliography}
\end{document}